\documentclass[12 pt]{amsart}

\usepackage{amssymb}

\usepackage[all]{xy}

\newtheorem{theorem}{Theorem}[section]
\newtheorem{lemma}[theorem]{Lemma}

\theoremstyle{definition}
\newtheorem{definition}[theorem]{Definition}
\newtheorem{example}[theorem]{Example}

\theoremstyle{remark}
\newtheorem{remark}[theorem]{Remark}

\numberwithin{equation}{section}

\begin{document}

\title[Slope of double cover fibrations]{A sharp bound for the slope of double cover fibrations}

\author{Maurizio Cornalba}
\address{Dipartimento di Matematica, Universit\`a di Pavia, Via Ferrata 1, 27100 Pavia, Italy}
\email{maurizio.cornalba@unipv.it}
\thanks{Research partially supported by: PRIN 2003 \textit{Spazi di moduli e teoria di Lie}; GNSAGA; FAR 2002 (Pavia) \textit{Variet\`a algebriche, calcolo algebrico, grafi orientati e topologici}}

\author{Lidia Stoppino}
\address{Dipartimento di Matematica, Universit\`a di Roma Tre, Largo San Leonardo Murialdo, 1, 00146 Roma, Italy}
\email{stoppino@mat.uniroma3.it}


\begin{abstract}
Let $f\colon X\rightarrow B$ be a fibration of genus $g$ whose general fiber is a double cover of a smooth curve of genus $\gamma$.
We show that $4(g-1)/(g-\gamma)$ is a sharp lower bound for the slope of $f$ when $g> 4\gamma+1$, proving a conjecture of Barja. 
Moreover, we give a 
characterization of the fibered surfaces that reach the bound. In the case $g=4\gamma+1 $ we obtain the same sharp bound under the additional assumption that the involutions on the general fibers glue to a global involution on $X$. 
\end{abstract}

\maketitle

\section*{Introduction and preliminaries}

A \emph{fibered surface}, or simply a \emph{fibration}, is a proper surjective morphism with connected fibers $f$ 
from a smooth surface $X$ to a smooth complete curve $B$. 
Call $F$ the general fiber of $f$.
A fibration is said to be \emph{relatively minimal} if the fibers contain no (-1)--curves.
The genus $g$ of $F$ is called the \emph{genus of the fibration}.
We say that $f$ is smooth if all the fibers are smooth, isotrivial if all the smooth fibers are mutually 
isomorphic, and \emph{locally trivial} if it is smooth and isotrivial.
A \emph{hyperelliptic} (respectively \emph{bielliptic}) \emph{fibration} 
is a fibered surface whose general fiber is a hyperelliptic (respectively bielliptic) curve.
A fibration is said to be \emph{semistable} if all its fibers are reduced nodal curves which are moduli semistable (any rational smooth component meets the rest of the curve in at least 2 points).

\subsubsection*{Relative invariants}

As usual, the \emph{relative dualizing
sheaf} of a fibration $f\colon X\rightarrow B$ is
the line bundle $$\omega_f =\omega_X\otimes (f^*\omega_B)^{-1},$$ where $\omega_V$ is the canonical
sheaf of $V$.

\begin{remark}
A relatively minimal fibration is a fibration such that $\omega_f$ is $f$-nef. A semistable fibration is a relatively minimal fibration whose fibers are nodal and reduced.
\end{remark}

The basic invariants for a relatively minimal fibration $f\colon X\rightarrow B$ are:
$$(\omega_f\cdot\omega_f)\,\mbox{ , } \,\,
\deg f_*\omega_f\,\,\mbox{ and }\,\,e_f:=e(X)-e(B)e(F),$$ 
where $e$ is the topological Euler number. They are related by Noether's formula:
$$
(\omega_f\cdot\omega_f)=12\deg f_*\omega_f-e_f.
$$
It is well known that all these invariants are greater or equal to $0$; moreover, 
$\deg f_*\omega_f=0$ if and only if $f$ is locally trivial, $(\omega_f\cdot\omega_f)=0$ implies that $f$ is
isotrivial, and $e_f=0$ if and only if the fibration is smooth.

\subsubsection*{Slope}

Assuming that the fibration is not locally trivial,
we can consider the ratio
$$
\mbox{s} (f)= \frac{(\omega_f\cdot\omega_f)}{\deg f_*\omega_f},
$$
which is itself an important invariant of the fibration, called the \emph{slope}. Noether's formula gives the upper bound $\mbox{s}(f)\leq 12$, which is achieved when all
the fibers are smooth (e.g., for the Kodaira fibrations). 
The lower bound for relatively minimal fibrations of genus $g\geq 2$ is given by the so-called
\emph{slope inequality} (see \cite{X}, \cite{C-H} and \cite{S})
\begin{equation}\label{slope}
(\omega_f\cdot\omega_f)\geq \frac{4(g-1)}{g} \deg f_*\omega_f;
\end{equation}
that is, $\mbox{s}(f)\geq 4(g-1)/g$.
Observe that the slope inequality implies in particular that $(\omega_f\cdot\omega_f)=0$ if and only if $f$ is locally trivial.
The inequality is sharp; equality holds for certain hyperelliptic fibrations (see \cite{C-H} and \cite{A-K}).

\medskip 

One of the main problems in the study of fibered surfaces is to understand how properties of the general fiber influence 
the slope. It is significant that, if the bound $4(g-1)/g$ is reached, then $F$ has a $g_2^1$. 
As a matter of fact, the gonality (or the Clifford index) of the 
general fiber plays an important role: there are several results in this direction 
(mainly due to Konno), although an explicit sharp bound on the slope depending on the gonality is still not 
known and seems to be hard to find.

\subsubsection*{Double covers}

Another direction in which the hyperelliptic case can be generalized is the study of fibrations whose general fiber $F$ 
is a double cover of a smooth curve of genus $\gamma\geq 0$, which we will call
double fibrations. 
As will be made clear in the first section, it is necessary to work with a smaller family of fibrations, i.e., fibrations that possess a \emph{global} involution which restricts to an involution of the general fibers (double cover fibrations).
The bielliptic case ($\gamma=1$) has been treated by Barja in 
\cite{B}: $4$ is the sharp lower bound for the slope and it is possible to give a characterization 
of the fibrations that reach the 
bound.
About the general case, what is known at present is a result of Barja and Zucconi, 
who show that the slope of a double cover fibration with  $g\geq 2\gamma+11$ is again greater than $4$ (cf. Theorem 0.6 in \cite{B-Z}).

It is easy to construct examples of double cover fibrations with slope $4(g-1)/(g-\gamma)$, 
while examples with smaller slope are known only for $g<4\gamma$ (see \cite{BPhD} and Examples \ref{primoesempio} and \ref{secondoesempio}). 
Moreover, for hyperelliptic and 
bielliptic fibrations the number $4(g-1)/(g-\gamma)$ gives exactly the sharp bound. 
It is therefore natural to conjecture, as Barja does in section 4.2 of
\cite{BPhD}, this to be the sharp bound for double cover fibrations with $g\geq 4\gamma+1$. 

In this paper, we give an affirmative answer to this conjecture (Theorem \ref{DC} and Theorem \ref{DC2}), together with a characterization
of the fibrations that reach the bound. Our results follow from an application of the slope inequality for fibered surfaces and of the algebraic index theorem, or signature theorem (see e.g. \cite{BHPVdV} Theorem IV.2.14 or \cite{G-H}).

In section \ref{dc} we discuss the problem of gluing involutions on the general fibers to a global involution. In section \ref{red} we recall a standard construction that allows to relate the invariants of a double cover fibration to the ones of a fibration which is a ``true'' double cover of a relatively minimal fibration of genus $\gamma$.
Section \ref{bound} is devoted
to the proof of the bound and to the characterization of the extremal case, while in section \ref{esempi} we present two examples  that prove the sharpness of the bound.

\section{Double covers and double fibrations}\label{dc}

The notion of double fibration is a natural extension of the one of hyperelliptic or bielliptic fibration.
\begin{definition} \label{df}
A \emph{double fibration of type $(g,\gamma)$} is a relatively minimal genus $g$ fibered surface $f\colon X\rightarrow B$ 
such that there is a degree $2$ morphism from the general fiber of $f$ to a smooth curve of genus $\gamma$.
\end{definition}

The sheet interchange involution on the general fiber of a double fibration $f\colon X\to B$ does not necessarily come from a global involution on $X$.
The problem, of course, is that the involution on the general fiber may not be unique, and there may not exist a rational section of $\mathbf{Aut}_{B}(X)\to B$, where $\mathbf{Aut}_{B}(X)$ stands for the relative Hilbert scheme parametrizing automorphisms of fibers of $X\to B$, reducing to the given involution on the general fiber. If such a section exists, for instance if the involution on the general fiber of $f$ is unique, it gives a global rational involution on $X$, which is actually regular when $g\ge 2$ because of the relative minimality of $X\to B$. In this case we thus get what we call a double cover fibration.

\begin{definition}\label{dcf}
A \emph{double cover fibration of type $(g,\gamma)$} is the datum of a genus $g$ fibration $f\colon X\rightarrow B$ together with a global involution on $X$ that restricts, on the general fiber, to an involution with genus $\gamma$ quotient. 
\end{definition}

This definition of double cover fibration is slightly more restrictive than the one given in \cite{B-Z}; in the cases we shall be concerned with, however, the two definitions are equivalent. As one can always produce sections of the scheme of relative automorphisms after a suitable finite base change $T\rightarrow B$, one might think that, to prove slope inequalities for general double fibrations of genus $g\ge 2$, it would suffice to prove them for double cover fibrations. While this strategy works fine if $X\to B$ is a semistable fibration, in the general case it runs into difficulties which at present appear insurmountable, as the slope behaves very badly under base change \cite{T1,T2}.

We are now going to show that, under the assumption $g> 4\gamma+1$, the involution on a general fiber of a double cover fibration of type $(g,\gamma)$ is indeed unique, and that the same happens when $g=4\gamma+1$, except in a very special case.
The argument is due to Barja (Lemma 4.7 in \cite{BPhD}), save for the discussion of the case $g=4\gamma+1$.

\begin{lemma}\label{invo} Let $F$ be a smooth curve of genus $g$, and let $\gamma\geq 1$ be an integer. If $g> 4\gamma+1$,
then $F$ has at most one involution $\iota$ such that $\Gamma=F/\langle\iota\rangle$ has genus $\le \gamma$. If instead $g=4\gamma+1$, and there are distinct involutions $\iota_1$ and $\iota_2$ of $F$ such that the quotients $F/\langle\iota_1\rangle=\Gamma_1$, $F/\langle\iota_2\rangle=\Gamma_2$ have genera $\gamma_1,\gamma_2$ both not exceeding $\gamma$, then $\gamma_1=\gamma_2=\gamma$, $\Gamma_1$ and $\Gamma_2$ are hyperelliptic, the natural map $F\to \Gamma_1\times\Gamma_2$ is an embedding, and its image belongs to the linear system $|\pi_1^*(2q_1)+\pi_2^*(2q_2)|$, where $q_i$ is a Weierstrass point on $\Gamma_i$ and $\pi_i$ denotes the projection $\Gamma_1\times\Gamma_2\to \Gamma_i$.
\end{lemma}
\begin{proof} Suppose $\iota_1$ and $\iota_2$ are two involutions of $F$ such that the quotients
$$
F/\langle\iota_1\rangle=\Gamma_1 ,\quad F/\langle\iota_2\rangle=\Gamma_2
$$
have genera $\gamma_1$ and $\gamma_2$ not greater than $\gamma$. Consider the commutative diagram
$$
\xymatrix@R-5pt{
 & F \ar[dl]_{\sigma_1} \ar[d]^\sigma \ar[dr]^{\sigma_2} & \\
 \Gamma_1 & D \ar[l]_{\beta_1} \ar[r]^{\beta_2} \ar@{^{(}->}[d]^j & \Gamma_2\\
 & \Gamma_1\times \Gamma_2 \ar[ul]^{\pi_1} \ar[ur]_{\pi_2} &
}
$$
where the $\sigma_i$ are the quotient morphisms, $D=\sigma_1\times\sigma_2(F)$, $j\circ\sigma=\sigma_1\times\sigma_2$, the $\pi_i$ are the projections, and the $\beta_i$ their restrictions to $D$. Clearly, the degree of $\sigma$ is either $1$ or $2$.
If it is $2$, the $\beta_i$ have to be isomorphisms; therefore $\sigma_1$ and $\sigma_2$ are the quotient maps of the \emph{same} involution on $F$. Conversely, if the involutions $\iota_1$ and $\iota_2$ coincide, the degree of $\sigma$ must be 2.

Now suppose that $\deg\sigma=1$.
Set $L_i=\pi_i^{-1}(p_i)\subseteq\Gamma_1\times\Gamma_2$, with $p_i\in\Gamma_i$. The effective divisor $L_1+L_2$ has self-intersection $2>0$.
By the index theorem  the determinant of the intersection matrix of the pair $(D,L_1+L_2)$ 
has to be non-positive. In other words
$$
2(D\cdot D)-(D\cdot L_1+L_2)^2\leq 0.
$$
Since $(D\cdot L_i)=2$, we obtain that $(D\cdot D)\leq 8$.
By the adjunction formula
\begin{equation}\label{adj}
2g-2 \leq 2p_a(D)-2 = (K_{\Gamma_1\times\Gamma_2}+D\cdot D)\leq
4(\gamma_1+\gamma_2)\leq 8\gamma.
\end{equation}
Thus $g\leq 4\gamma +1$, and the first part of the lemma is proven.

If $g=4\gamma +1$, all the above inequalities must necessarily be equalities. 
In particular, $\gamma_1=\gamma_2=\gamma$, and $p_a(D)=g$, so $\sigma$ is an isomorphism.
Furthermore, by the index theorem, $D$ has to be numerically equivalent to a rational multiple of $L_1+L_2$. 
Intersecting with $L_1$ and $L_2$, one sees that $D$ is numerically equivalent to $2L_1+2L_2$.
Hence $D$ is linearly equivalent to a divisor of the form $\pi_1^*A_1+\pi_2^*A_2$, where $A_i$ is a divisor of degree $2$ on $\Gamma_i$.
Observe that 
\begin{equation}\label{kunneth}
H^0(\Gamma_1\times\Gamma_2\, ,\mathcal O(\pi_1^*A_1+\pi_2^*A_2))=
H^0(\Gamma_1,\mathcal O(A_1))\otimes H^0(\Gamma_2,\mathcal O(A_2)).
\end{equation}
It follows in particular that the dimension of $H^0(\Gamma_i,\mathcal O(A_i))$ must be strictly positive, since the linear system $|\pi_1^*A_1+\pi_2^*A_2|$ is non-empty. If the dimension of $H^0(\Gamma_1,\mathcal O(A_1))$ were equal to $1$, formula \ref{kunneth} would imply that every divisor in $|\pi_1^*A_1+\pi_2^*A_2|$ is of the form $\pi_1^*D_1+\pi_2^*D_2$, where $D_1$ is the unique divisor in $|A_1|$ and $D_2$ belongs to $|A_2|$, and is therefore singular. Since $D$ is smooth and belongs to $|\pi_1^*A_1+\pi_2^*A_2|$, it follows that the dimension of $H^0(\Gamma_1,\mathcal O(A_1))$ must be at least $2$. Similar considerations apply to $H^0(\Gamma_2,\mathcal O(A_2))$. We therefore conclude that the $\Gamma_i$ are hyperelliptic and that $A_i$ is linearly equivalent to $2 q_i$, where $q_i$ is a Weierstrass point.
\end{proof}

Lemma \ref{invo} implies in particular that, if $g=4\gamma+1$ and $F$ has an involution $\iota$ such that $F/\langle\iota\rangle$ has genus $\gamma$, the involution is unique if $F/\langle\iota\rangle$ is non-hyperelliptic.

If the general fibers of a double fibration are as described in the second part of Lemma \ref{invo}, it is possible that a non-trivial base change is needed in order to get a global involution, as the following example shows (see also \cite{B} for an example in the bielliptic case).

\begin{example}
Let $\Gamma$ be a smooth hyperelliptic curve of genus $\gamma$. Let $B$ be a smooth curve of positive genus, and let $\alpha\colon T\to B$ be an unramified degree two covering; thus $B$ is the quotient of $T$ modulo a base-point-free involution $\sigma$. We set $G=\langle\sigma\rangle$, and denote by $\tau$ the exchange of components automorphism of $\Gamma\times \Gamma$. Call $\pi_1$, $\pi_2$ the projections from $\Gamma\times \Gamma$ to the two factors. Let $q$ be a Weierstrass point of $\Gamma$, and consider  the linear system $|\pi_1^*(2q)+\pi_2^*(2q)|^\tau$ of effective $\tau$-invariant divisors linearly equivalent to $\pi_1^*(2q)+\pi_2^*(2q)$; it is immediate to show that it is  base-point-free and not composed with an involution. So, by Bertini's theorem, a general member $F\in |\pi_1^*(2q)+\pi_2^*(2q)|^\tau$ is smooth and irreducible. The genus of $F$ is $g=4\gamma+1$.
Let $G$ act on $\Gamma\times \Gamma\times T$ by
$$
\sigma (f_1,f_2,t)=(f_2, f_1, \sigma(t))=(\tau(f_1,f_2),\sigma(t)).
$$
Clearly, the subvariety $F\times T\subseteq \Gamma\times \Gamma\times T$ is $G$-invariant, and the action of $G$ on it is compatible with the action of $G$ on $T$. Dividing by these two actions, we thus get from $F\times T\to T$ a new fibration $X\to B$, where $X=(F\times T)/G$, with fibers all isomorphic to $F$.  The curve $F$ carries two involutions $\iota_1,\iota_2$, corresponding to the projections to the two factors of $\Gamma\times\Gamma$. Notice that $\iota_1$ and $\iota_2$ are conjugate under the involution $\tau$, i.e., that $\tau\iota_1=\iota_2\tau$. For each $b\in B$, the fiber $X_b$ of $X\to B$ inherits from $F$ two involutions with quotient of genus $\gamma$. In particular, $X\to B$ is a double fibration of type $(g,\gamma)$. We claim that, for any $b\in B$, the two involutions on $X_b$ belong to the same component of the scheme $\mathbf{Aut}_{B}(X)$. In fact, $\mathbf{Aut}_{B}(X)$ is the quotient modulo $G$ of $\mathbf{Aut}_{T}(F\times T)={\rm Aut}(F)\times T$, where $\sigma$ acts by sending $(\alpha,t)\in {\rm Aut}(F)\times T$ to $(\tau\alpha\tau,\sigma(t))$, while the irreducible components of $\mathbf{Aut}_{B}(X)$ are the images of the irreducible components of $\mathbf{Aut}_{T}(F\times T)$. On the other hand, the images of the components $\{\iota_1\}\times T$ and $\{\iota_2\}\times T$ coincide, as follows from the observation that $\sigma(\iota_1,t)=(\tau\iota_1\tau,\sigma(t))=(\iota_2,\sigma(t))$. This proves our claim. A consequence is that, for any $b\in B$, there is no global involution $\iota$ in the group ${\rm Aut}_B(X)$ of automorphisms of $X$ over $B$ which pulls back on $X_b$ to one of the two involutions coming from $\iota_i$ and $\iota_2$. In fact, if such a $\iota$ existed, it would give an involution on every fiber of $X\to B$, and hence a section of $\mathbf{Aut}_{B}(X)$ over $B$. This section would be a component of $\mathbf{Aut}_{B}(X)$ containing only one of the two involutions on $X_b$, contrary to what we just proved.

\end{example}

\section{Reduction to double covers of relatively minimal fibrations}\label{red}

As this is crucial for our argument, we now sketch a well-known procedure that 
associates to a double cover fibration a double cover of a  relatively minimal fibration. 
The precise construction can be found for instance in \cite{A-K} or in \cite{B-Z}.

Let $f\colon X\rightarrow B$ be a double cover fibration. Let $\iota$ be the involution on $X$. 
If it has a fixed locus of codimension $1$, the quotient 
$X/\langle\iota\rangle$ is a smooth surface. 
Otherwise, consider the blow-up $\widetilde X$ of $X$ at the isolated fixed points of $\iota$,
and call $\tilde\iota$ the induced involution on it.
The quotient $\widetilde X/\langle \tilde\iota\rangle=\widetilde Y$ is a smooth surface with a natural fibration $\widetilde\alpha$
over $B$ which is not necessarily relatively minimal. Let $\alpha\colon Y\to B$ be its minimal model.
$$
\xymatrix@R-7pt{
\widetilde X \ar[d] \ar[r] & \widetilde Y \ar[d]_{\widetilde\alpha} \ar[r] &Y \ar[dl]^\alpha\\
X \ar[r]^f & B &\\
}
$$
The direct image  $R$ of the branch locus of the double cover 
$\widetilde X\rightarrow \widetilde Y$
induces a double cover $X'\rightarrow Y$, with $X'$ normal but not
necessarily smooth; notice however that, by construction, $X'$ is a hypersurface in a threefold which is smooth over $B$, so $X'\to B$ admits an invertible relative dualizing sheaf.
To obtain a smooth double cover we perform the {\em canonical resolution}
(see \cite{BHPVdV} III.7, \cite{B} sec. 2, \cite{A-K} sec. 2.2)
$$
\xymatrix@R-7pt{
 X_k \ar[r]^{\sigma_k} \ar[d] & X_{k-1} \ar[r] \ar[d] & \cdots & X_1  \ar[r]^{\sigma_1} \ar[d] & X_0 =X' \ar@<-12pt>[d] \\
Y_k \ar[r]^{\tau_k} & Y_{k-1} \ar[r] & \cdots & Y_1 \ar[r]^{\tau_1} & Y_0 =Y\\
}
$$ 
where the $\tau_j$ are successive blow-ups that resolve the singularities of $R$; 
the morphism $X_j\rightarrow Y_j$
is the double cover with branch locus 
$R_j:=\tau_j^*R_{j-1}-2\left[\frac{m_{j-1}}{2}\right]E_j$, where $E_j$ is the exceptional 
divisor of $\tau_j$, $m_{j-1}$ is the multiplicity of the blown-up point, and $[\ \ ]$ stands for integral part.
Let $f_j\colon X_j\rightarrow B$, $f'\colon X'\rightarrow B$  be the induced fibrations. 
A computation shows that
$$
(\omega_{f_k}\cdot \omega_{f_k})
=(\omega _{f'}\cdot\omega _{f'})-2\sum_{i=1}^{k}\left(\left[\frac{m_i}{2}\right]-1\right)^2,
$$
and that
$$
\deg (f_k)_*\omega_{f_k}=
\deg f'_*\omega _{f'}-\frac{1}{2}\sum_{i=1}^{k}\left[\frac{m_i}{2}\right]\left(\left[\frac{m_i}{2}\right]-1\right).
$$
Observe that, since $X_k$ is smooth, 
by the relative minimality of $f\colon X\rightarrow B$ there is a morphism
$\beta \colon X_k\rightarrow X$. Therefore
$$
(\omega_f\cdot\omega_f)=(\omega_{f_k}\cdot\omega_{f_k})+ \epsilon,
$$
where $\epsilon$ is the number of blow-ups which make up $\beta$.
Moreover, observe that $f_*\omega_f=(f_k)_*\omega_{f_k}$.
Hence we get the following fundamental identity:
\begin{equation}\label{reduction}
\begin{aligned}
&(\omega_f\cdot\omega_f)-4\frac{g-1}{g-\gamma}f_*\omega_f \\
&=(\omega_{f'}\cdot\omega_{f'})-4\frac{g-1}{g-\gamma}f'_*\omega_{f'}
+2\sum_{i=1}^{k}
\left(\left[\frac{m_i}{2}\right]-1\right)\left(\frac{\gamma-1}{g-\gamma}\left[\frac{m_i}{2}\right]+1\right)+\epsilon.
\end{aligned}
\end{equation}

\begin{definition}
In the situation above, we say that the branch divisor $R\subset Y$ has \emph{negligeable singularities} if all the multiplicities $m_i$ in the above process equal 2 or 3 (cf. \cite{P})
\end{definition}

\section{The bound}\label{bound}

\begin{theorem}\label{DC}
Let $f\colon X\rightarrow B$ be a double fibration of type $(g,\gamma)$. If $g> 4\gamma+1$, then
\begin{equation}\label{ine}
(\omega_f\cdot \omega_f)\geq 4\frac{g-1}{g-\gamma}\deg f_*\omega_f.
\end{equation}
If $\gamma\geq 1$, equality holds if and only if $X$ is the minimal
desingularization of a double cover $\pi\colon  \overline X\rightarrow Y$
of a locally trivial genus $\gamma$ fibration $\alpha\colon  Y\rightarrow B$ 
such that the branch locus $R$ of $\pi$ has only negligeable singularities and, in addition, when $\gamma>1$, is numerically equivalent to a rational linear combination of $\omega_\alpha$ and a fiber of $\alpha$.
\end{theorem}

\begin{proof} The case $\gamma=0$ is the slope inequality for hyperelliptic fibrations.
The case $\gamma=1$ has been proved in \cite{B}.
We therefore assume that $\gamma>1$.
In view of Lemma \ref{invo}, the assumptions about $g$ and $\gamma$ guarantee that we are in fact dealing with a double cover fibration. 
We adopt the notation introduced in the previous section.
In view of the identity (\ref{reduction}), to prove (\ref{ine}) it suffices to prove its analogue for $f'\colon X'\to B$. Recall that the double covering $\xi\colon X'\to Y$ corresponds to a line bundle $\mathcal L$ on $Y$ such that $\mathcal L^2=\mathcal O(R)$, where $R$ is the ramification divisor of $\xi$, and that
$$
\xi_*\mathcal O_{X'}=\mathcal O_Y\oplus \mathcal L^{-1}\,,\quad
\omega_{f'}=\xi^*(\omega_\alpha\otimes\mathcal L)\,.
$$ 
It follows that
$$
(\omega_{f'}\cdot \omega_{f'})=2(\omega_\alpha\otimes\mathcal L\cdot \omega_\alpha\otimes\mathcal L)=2(\omega_\alpha\cdot \omega_\alpha)+4(\mathcal L\cdot \omega_\alpha)+2(\mathcal L\cdot \mathcal L),
$$
and also, by the Riemann-Roch theorem, that
$$
 \deg {f'}_*\omega_{f'} =2\deg\alpha_*\omega_\alpha+\frac{(\mathcal L\cdot \mathcal L)}{2}+\frac{(\mathcal L\cdot \omega_\alpha)}{2}.
$$
Hence we may write:
$$
\begin{aligned}
&(\omega_{f'}\cdot \omega_{f'})-4\frac{g-1}{g-\gamma}\deg {f'}_*\omega_{f'}\\
&=2\left((\omega_\alpha\cdot \omega_\alpha)-4\frac{g-1}{g-\gamma}\deg\alpha_*\omega_\alpha\right) 
- 2\frac{\gamma-1}{g-\gamma}(\mathcal L\cdot \mathcal L)+2\frac{g-2\gamma+1}{g-\gamma}(\mathcal L\cdot \omega_\alpha).
\end{aligned}
$$
Using the slope inequality (\ref{slope}) for $\alpha\colon  Y\rightarrow B$ we obtain that
$$
(\omega_\alpha\cdot \omega_\alpha)-4\frac{g-1}{g-\gamma}\deg\alpha_*\omega_\alpha\geq \frac{-g-\gamma^2+2\gamma}{(\gamma-1)(g-\gamma)}(\omega_\alpha\cdot \omega_\alpha).
$$
Therefore
$$
\begin{aligned}
&(\omega_{f'}\cdot \omega_{f'})-4\frac{g-1}{g-\gamma}\deg {f'}_*\omega_{f'}\\
&\geq\frac{2}{g-\gamma}\left( (g-2\gamma+1)(\omega_\alpha\cdot \mathcal L)-(\gamma-1)(\mathcal L\cdot \mathcal L)-\frac{\gamma^2+g-2\gamma}{\gamma-1}(\omega_\alpha\cdot \omega_\alpha)\right)\\
&=\frac{1}{g-\gamma}\left( 2(\mathcal L\cdot \Gamma)(\omega_\alpha\cdot \mathcal L)-(\omega_\alpha\cdot \Gamma)(\mathcal L\cdot \mathcal L)-4\frac{\gamma^2+g-2\gamma}{(\omega_\alpha\cdot \Gamma)}(\omega_\alpha\cdot \omega_\alpha)\right),
\end{aligned}
$$
where $\Gamma$ stands for a general fiber of $\alpha$. 
As $(\omega_\alpha\cdot\omega_\alpha)\geq 0$, the intersection matrix of $\omega_\alpha$, $\mathcal L$ and $\Gamma$ 
cannot be negative definite. 
The index theorem then implies that its determinant is non-negative, i.e., that
$$
2(\mathcal L\cdot \Gamma)(\omega_\alpha\cdot \Gamma)(\omega_\alpha\cdot \mathcal L)-(\omega_\alpha\cdot \Gamma)^2(\mathcal L\cdot \mathcal L)\geq (\mathcal L\cdot \Gamma)^2(\omega_\alpha\cdot \omega_\alpha).
$$
Combining this inequality with the ones obtained above, we get
$$
(\omega_{f'}\cdot \omega_{f'})-4\frac{g-1}{g-\gamma}\deg {f'}_*\omega_{f'}\geq
\frac{1}{g-\gamma}\left(\frac{(\mathcal L\cdot \Gamma)^2}{(\omega_\alpha\cdot \Gamma)}-4\frac{\gamma^2+g-2\gamma}{(\omega_\alpha\cdot \Gamma)}\right)(\omega_\alpha\cdot \omega_\alpha),
$$
and so
\begin{equation}\label{fond}
(\omega_{f'}\cdot \omega_{f'})-4\frac{g-1}{g-\gamma}\deg {f'}_*\omega_{f'}\geq\frac{(g-4\gamma-1)(g-1)}{2(g-\gamma)(\gamma-1)}(\omega_\alpha\cdot \omega_\alpha).
\end{equation}
The expression on the right is clearly non-negative as soon as $g\ge 4\gamma+1$.
Note that the argument so far applies to any double cover fibration.

To prove the characterization of the fibrations that reach the bound, observe first that the coefficient of $(\omega_\alpha\cdot \omega_\alpha)$ in (\ref{fond}) is not zero when $g>4\gamma+1$, so the local triviality of $\alpha$ is a necessary condition. 
Then recall that $\mathcal O(R)=\mathcal L^2$ and notice that, if (\ref{ine}) is an equality, all the inequalities in the proof must be equalities, and the terms $2\sum_{i=1}^{k}
\left(\left[\frac{m_i}{2}\right]-1\right)\left(\frac{\gamma-1}{g-\gamma}\left[\frac{m_i}{2}\right]+1\right)$ and $\epsilon$ in (\ref{reduction}) must vanish. In particular, we get that
$$
(g-2\gamma+1)(\omega_\alpha\cdot \mathcal L)-(\gamma-1)(\mathcal L\cdot \mathcal L)=0
$$
which, in view of of $(\omega_\alpha\cdot\omega_\alpha)=0$ and of $\mathcal O(R)=\mathcal L^2$, is equivalent to the vanishing of the determinant of the intersection matrix of $\omega_\alpha$, $\Gamma$ and $R$, i.e., to $R$ being numerically equivalent to a rational linear combination of $\omega_\alpha$ and $\Gamma$.
\end{proof}

The analogous result for $g=4\gamma +1$ can be stated as follows.

\begin{theorem}\label{DC2}
Let $f\colon X\rightarrow B$ be a double fibration of type $(g,\gamma)$ with $g=4\gamma +1$. Then inequality (\ref{ine}) holds, provided we are in one of the following cases:
\begin{enumerate}
\item $f$ is a double cover fibration;
\item $f$ is a semistable fibration.
\end{enumerate}
In particular, (\ref{ine}) is valid if a smooth fiber of $f$ admits an involution whose quotient is a non-hyperelliptic curve of genus $\gamma$.

Moreover, a necessary condition for the slope to reach the bound is that the associated relatively minimal fibration of genus $\gamma$ be either locally trivial or hyperelliptic with slope $4(\gamma-1)/\gamma$.
\end{theorem}
\begin{proof} Case (1) is covered by the argument used to prove Theorem \ref{DC}.
In case (2), let $\rho\colon T\to B$ be a finite map with $T$ smooth. Then $X\times_BT$ has $A_n$ singularities. Let $X'$ be its minimal resolution, and $f'\colon X'\to T$ the natural projection; clearly, this is a minimal fibration. As is well known, the slope of $f'$ is equal to the one of $f$. In fact, in this case, $\omega_{f'}$ is just the pullback to $X'$ of $\omega_f$, so $(\omega_{f'}\cdot \omega_{f'})$ and $\deg {f'}_*\omega_{f'}$ are both equal to $\deg\rho$ times the corresponding invariant of $f$. Since the base change $\rho$ can be chosen so that $f'\colon X'\to T$ is a double cover fibration, we are reduced to the previous case. It follows from Lemma \ref{invo} and the comment immediately following its proof that a sufficient condition for $f$ to be a double cover fibration is that \emph{one} of its smooth fibers admit an involution with non-hyperelliptic genus $\gamma$ quotient.
The coefficient of $(\omega_\alpha\cdot\omega_\alpha)$ in inequality (\ref{fond}) is $0$ when $g=4\gamma +1$. 
Hence the local triviality of $\alpha$ is no longer a necessary condition for the fibration to reach the bound.
If $\alpha$ is not locally trivial, it is instead necessary that $\alpha$ itself attain the bound given by the slope inequality, so we conclude.
\end{proof}

Clearly, one could give  necessary and sufficient conditions as in Theorem \ref{DC}, imposing that the inequalities in the proof be equalities.
It is interesting to notice that, in this borderline case, the conditions change substantially, because local triviality of the fibration of genus $\gamma$ is no more needed, and indeed one can construct a fibered surface of arbitrary genus $g$ reaching the bound and which is a double cover of a non locally trivial fibration of genus $\gamma=(g-1)/4$ (Example \ref{secondoesempio}).

\section{Examples}\label{esempi}

We present below two examples, both due to Barja (cf. \cite{BPhD}, sec. 4.5) which show that the bound given is indeed sharp. 
The first is an example of double cover fibration reaching the bound; in the second we construct a fibration with $g=4\gamma+1$, reaching the bound, which is a double cover of a hyperelliptic fibration which in turn reaches the bound given by the slope inequality. 
The last  construction also leads to counterexamples to the bound for $g<4\gamma$.

\begin{example}\label{primoesempio}
This is a generalization of the examples constructed in \cite{X} and in \cite{C-H}. Let $\Gamma$ and $B$ be smooth curves. 
Call $\gamma$ the genus of $\Gamma$.
Let $p_1\colon B\times \Gamma\rightarrow B$ and $p_2\colon B\times \Gamma\rightarrow \Gamma$ be the two projections, and $H_1$, $H_2$ their general fibers.
For sufficiently large integers $n$ and $m$, the linear system $|2nH_1+2mH_2|$ is base-point-free. 
Hence, by Bertini's Theorem there exists a smooth divisor $R\in |2nH_1+2mH_2|$. 
As $R$ is even, we can construct the double cover $\rho\colon X\rightarrow B\times\Gamma$ ramified 
over $R$. 
Consider the fibration $f:=p_1\circ\rho \colon X\rightarrow B$; 
its  general fiber is a double cover of $\Gamma$, and its genus is $g=2\gamma+m-1$.
Observe that
$$
\omega_f\cong\rho^*(\omega_{p_1}(nH_1+mH_2))\cong \rho^*(\mathcal O(nH_1+(2\gamma-2+m)H_2)),
$$
and
$$
\deg f_*\omega_f =\deg {p_1}_*(\omega_{p_1}(nH_1+mH_2))=n(\gamma-1+m).
$$
Therefore the slope of $f$ is exactly
$$ s(f)=4\frac{2\gamma+m-2}{\gamma+m-1}=4\frac{g-1}{g-\gamma}.$$ 
If we consider a general divisor $R\in |2nH_1+2mH_2|$, it
has only simple ramification points over $B$, and we obtain a semistable fibration.
Notice that $R$ is numerically equivalent to a linear combination of $\omega_{p_1}$ and $\Gamma$, as it should be, because
$$
R\equiv 2nH_1+2mH_2\equiv 2n\Gamma+\frac{m}{\gamma-1}K_{p_1}.
$$
\end{example}

\begin{example}\label{secondoesempio}
Consider the fibration of Example \ref{primoesempio} with $\Gamma=\mathbb P^1$, and set $f_i=p_i\circ\rho$.
Call $F_i$ the general fiber of $f_i$; hence $F_1$ is hyperelliptic of genus $\gamma=m-1$. 
By what we observed in Example \ref{primoesempio},
$$
\omega_{f_1}\cong \rho^*(\omega_{p_1}(nH_1+mH_2))\cong \mathcal O(nF_1+(m-2)F_2).
$$ 
Let $x,y$ be positive integers, and consider the linear system $|2xF_1+2yF_2|$. 
Applying Bertini's theorem again, for large enough $x$ and $y$ we can find a smooth even divisor $\Delta$ belonging to it.
Let $\pi \colon Y\rightarrow X$ be the double cover ramified over $\Delta$.
Call $h$ the fibration $p_1\circ\rho\circ\pi \colon Y\rightarrow B$; the general fiber $F$ of $h$ is a double cover of $F_1$. Its genus is $g=2(m-1)+2y-1$.
Now, $\omega_h\cong\pi^*(\omega_{f_1}(xF_1+yF_2))\cong\mathcal O((n+x)F+(m+y-2)\pi^*F_2)$, so
$$
(\omega_h\cdot\omega_h)=8(x+n)(y+m-2),
$$
while
$$
\begin{aligned}
h_*\omega_h&={f_1}_*(\omega_{f_1}(xF_1+yF_2))\oplus {f_1}_*(\omega_{f_1})\\
&={p_1}_*(\rho_* (\rho^* \mathcal O((n+x)H_1+ (m+y-2)H_2)))\\
&={p_1}_*(\omega_{p_1}((n+x)H_1+(m+y)H_2))\oplus {p_1}_*\mathcal O(xH_1+(y-2)H_2)\oplus {f_1}_*(\omega_{f_1}).
\end{aligned}
$$
Therefore $$\deg h_*\omega_h=(x+n)(y+m-1)+x(y-1)+n(m-1).$$
If we choose, as we may, $m=y$, we get exactly $g=4m-3=4\gamma+1$ and slope
$$
s(h)=8\frac{2m-2}{3m-2}=4\frac{g-1}{g-\gamma}.
$$
Notice moreover that choosing $m>y$ we obtain fibrations with $g\leq 4\gamma-1$ and slope strictly smaller than $4(g-1)/(g-\gamma)$.
\end{example}

\bibliographystyle{amsplain}
\providecommand{\bysame}{\leavevmode\hbox to3em{\hrulefill}\thinspace}

\end{document}